\documentclass[letterpaper,12pt]{article}
\title{\bf Solving the Likelihood Equations}

\author{Serkan Ho\c sten, Amit Khetan  and Bernd Sturmfels}

\date{}

\usepackage{latexsym,array,delarray,amsthm,amssymb,epsfig}
\usepackage[sumlimits]{amsmath}

\theoremstyle{plain}
\newtheorem{thm}{Theorem}
\newtheorem{alg}[thm]{Algorithm}

\newtheorem{prop}[thm]{Proposition}

\theoremstyle{definition}
\newtheorem{df}[thm]{Definition}
\newtheorem{ex}[thm]{Example}
\newtheorem{rem}[thm]{Remark}
\newtheorem{pr}[thm]{Problem}

\def\eea{\end{eqnarray*}}
\def\bea{\begin{eqnarray*}}

\newcommand{\nn}{\mathbb{N}}
\newcommand{\pp}{\mathbb{P}}
\newcommand{\qq}{\mathbb{Q}}
\newcommand{\rr}{\mathbb{R}}
\newcommand{\cc}{\mathbb{C}}

\begin{document}

\maketitle

\begin{abstract} 
Given a model in algebraic statistics and data, 
the likelihood function is a rational function
on a projective variety.  Algebraic algorithms are presented 
for computing all critical points of this function,
with the aim of identifying the local maxima in the
probability simplex. 
Applications include models specified by rank conditions on matrices
and the Jukes-Cantor models of phylogenetics.
The maximum likelihood degree of a generic
complete intersection is also determined.
\end{abstract}

\section{Introduction}

A model in algebraic statistics is specified by a polynomial 
map from the space of model parameters to the space
of the joint probability distributions of the observed 
discrete random variables.
{\em Maximum likelihood estimation} is concerned with
 finding those model parameters that best explain
 a given sequence of observations. This is done 
 by maximizing the likelihood function.
The likelihood function is usually not convex,
it can have many local maxima, and the problem of
finding and certifying a global maximum is difficult.

Here we consider the problem of finding {\bf all}
critical points of the likelihood function,
with the aim of identifying all local maxima.
The defining equations of the
critical points are the {\em likelihood equations}.
The number of {\bf complex} solutions 
to the likelihood equations (for generic data) is  called the  
{\em maximum likelihood (ML) degree} 
of the model. A geometric study of the ML degree
was undertaken in  our joint work
with Fabrizio Catanese \cite{CHKS}. The present paper offers
algebraic algorithms for deriving and solving the likelihood equations.

We begin by illustrating the problem and our solution
for a simple example. In a certain game of chance,
a gambler tosses the same  coin four times in a row, and
the number of times heads come up are recorded. Hence the possible 
outcomes are $0$, $1$, $2$, $3$, or $4$.
 We observe $1000$ rounds of this game,
and we record the outcomes in the data vector
$ u \, = \, (u_0,u_1,u_2,u_3,u_4) \,\in\, \nn^5 $,
where $u_i$ is the number of trials that had $i$ heads.
Hence $u_0+u_1+u_2+u_3+u_4 = 1000$. Suppose we
are led to suspect that the gambler uses two biased coins, 
one in each of his sleeves, and he picks the coin to be used
at random (with probabilities $\pi$ and $1-\pi$)
prior to each round.
We wish to test this hypothesis using the data~$u$.

Our model is  the mixture of a pair of four-times repeated Bernoulli trials.
The mixing parameter $\pi$ is the probability that the gambler picks
the coin in his left sleeve. The bias of the left coin is $s$, and the
bias of the right coin is $t$. Our model stipulates that 
the probabilities of the five outcomes are
 \begin{eqnarray*}
p_0 \,\, = & \pi  (1-s)^4 \, \, + \,\, (1-\pi)  (1-t)^4  ,\\
p_1 \,\, = & 4 \pi  s (1-s)^3 \, \, + \,\, 4 (1-\pi)  t (1-t)^3,  \\
p_2 \,\, = & 6 \pi  s^2 (1-s)^2 \, \, + \,\, 6 (1-\pi)  t^2 (1-t)^2 , \\
p_3 \,\, = & 4 \pi  s^3 (1-s) \, \, + \,\, 4 (1-\pi)  t^3 (1-t) ,\\
p_4 \,\, = & \pi  s^4 \, \, + \,\,  (1-\pi)  t^4.
\end{eqnarray*}
 The polynomial $p_i$ represents the probability of seeing $i$ heads in 
a round.
The likelihood of observing the data $u$
 when $1000$ trials are made equals
\begin{equation}
\label{likelihoodfct}
\frac{ 1000 !  }{ u_0 ! u_1 ! u_2 ! u_3 ! u_4 !} \cdot
 p_0^{u_0} p_1^{u_1} p_2^{u_2} p_3^{u_3} p_4^{u_4} 
\end{equation}
Maximum likelihood estimation means
maximizing (\ref{likelihoodfct}) subject to  $0 < \pi,s,t < 1$.
The critical equations for this {\em unconstrained optimization problem}
have infinitely many solutions: there is a curve
of critical points in the $s=t$ plane. 

In order to avoid such non-identifiability, we
reformulate our maximum likelihood computation as 
the following {\em constrained optimization problem}:
\begin{equation}
\label{constrainedprob}
 {\rm Maximize} \,\,
  p_0^{u_0} p_1^{u_1} p_2^{u_2} p_3^{u_3} p_4^{u_4}
\,\,  \hbox{subject to} \,\,
  {\rm det}({\bf P}) = 0 \,\, \hbox{and} \,\, 
 p_0+ \cdots + p_4 = 1, \,\,\,\,
 \end{equation}
$$   {\rm where } \qquad {\bf P} \,\,\, = \,\,\,
\bmatrix
12 p_0 & 3  p_1 & 2 p_2 \\
 3 p_1 & 2 p_2 & 3 p_3 \\
 2 p_2 & 3 p_3 & 12 p_4
\endbmatrix.
$$
The image of the map $(\pi,s,t) \mapsto (p_0,p_1,p_2,p_3,p_4)$
over the complex numbers
is the  hypersurface $\{{\rm det}({\bf P})=0\}$
in projective $4$-space.
Using Algorithm \ref{implicit}, we find that
the ML degree of this model is $12$, i.e., the
solution of problem (\ref{constrainedprob})
leads to  an algebraic equation of degree $12$. 
See Examples~\ref{mixBern} and \ref{Bernoulli1}.

This paper is organized as follows. In Section 2 we introduce 
the likelihood equations associated with an arbitrary projective variety $V$.
The ML degree of $V$ is defined as the
number of complex solutions to the likelihood equations.

Section 3 contains an algebraic geometry result.
An explicit formula is given  for the ML degree of a generic
complete intersection. This formula is an upper bound
for the ML degree of more special complete intersections.

In Section 4 we present an algorithm whose input is
an arbitrary homogeneous ideal in a polynomial ring, 
representing a projective variety $V$. The algorithm
uses linear algebra over the coordinate ring $\rr[V]$
to find the  likelihood ideal.
This ideal typically has finitely many complex solutions.
We also discuss our test implementation in {\sc Singular} \cite{GPS01}.
It computes all solutions numerically and identifies
the local maxima in the probability simplex.

Section 5 comprises an experimental study of the
ML degree and number of local maxima 
for various determinantal models, including
the one discussed above. It is important to note that, in the 
context of algebraic statistics, every variety $V$ comes with a 
fixed coordinate system.  We demonstrate that the ML degree is 
extremely sensitive to changes of coordinates, even just scaling of the
coordinates. The good news is that in each case the ML degree
appears to be smallest for the statistically meaningful
coordinate system.

In Section 6 we apply our results to a class of models widely used in
computational biology: Jukes-Cantor models for
phylogenetic trees \cite{CHS, CHPH, CKS}.

 The setup of Sections 2--4 assumes that
 the defining ideal of the model $V$ is known. If this ideal is not known
and impossible to compute, then we are confined to use the
(generally less efficient)
parametric version of the likelihood equations which are
discussed in Section 7. In that section we also prove
that the parametric ML degree 
(which is the quantity emphasized in \cite{CHKS}) equals the implicit
ML degree times the cardinality of a generic fiber.

\section{Likelihood Locus on a Projective Variety}

We consider a statistical model which is a subset
of the probability simplex
$$ \Delta_n \,\,\,\, = \,\,\,\, \bigl\{ (p_0,p_1,\ldots,p_n) \in \rr^{n+1} \, :\,\,
p_0,\ldots,p_n > 0 \,\, \hbox{and} \, \,
p_0+p_1+\cdots + p_n \, = 1 \bigr\}, $$
and we assume that the model is presented as the solution set in $\Delta_n$
of a system of homogeneous polynomial equations in the unknowns $p_0,p_1,\ldots,p_n$.
Such polynomials are known as {\em model invariants} in the literature on phylogenetics 
and algebraic statistics \cite{StSu}. We write $V$ for the Zariski closure of the model
in complex projective space $\pp^n$. Equivalently, $V$ is the set of all complex solutions
to the given homogeneous polynomial equations. 
The  {\em maximum likelihood problem} 
is to find a point $p = (p_0:\cdots:p_n)$ in the model 
$$ V_{>  0} \,\,\, = \,\,\, V \,\,\cap \,\, \Delta_n $$
which ``best explains'' a given data vector 
$(u_0, \ldots, u_n) \in \nn^{n+1}$. As in
(\ref{constrainedprob}) above,
this means solving the following constrained optimization problem:
\begin{equation}
\label{mleproblem} {\rm Maximize} \,\, \,L \,= \, \frac{
p_0^{u_0} p_1^{u_1}  \cdots p_n^{u_n} }{ (p_0+p_1+\cdots +
p_n)^{u_0+u_1+\cdots + u_n}} \quad \hbox{subject to} \,\,
p \in V_{> 0}. 
\end{equation}

Our approach is to compute all complex 
critical points of the {\em likelihood function} $L$ and to
extract the positive real solutions that are local maxima.  
While the optimization problem (\ref{mleproblem})
requires the $p_i$ to be real and positive, 
we shall compute all the  critical  points on the 
complex projective  variety $V$.
Let $V_{sing}$ denote the singular locus of the variety $V$
and set $\,V_{reg} := V \backslash V_{sing}$. 
Let $P$ be the homogeneous ideal in the polynomial ring
$\rr[p_0,p_1,\ldots,p_n]$ generated by
the defining polynomials of $V$. All computations in the coordinate ring
$$ \rr[V] \quad = \quad \rr[p_0,p_1,\ldots,p_n]/P $$
will be made using standard techniques of Gr\"obner basis theory \cite{CLO, GP}.

\begin{df} Let $\mathcal{U} $ be the open subset 
$\,V_{reg} \setminus \mathcal{V}(p_0\cdots p_n \cdot (\sum p_i))\,$
of $\, V \subset \pp^n$. The {\em likelihood locus} $Z_u$ is
the set of points $p \in
\mathcal{U}$ such that $dL = 0$. The {\em likelihood ideal} $I_u
\subset \rr[V]$ is the ideal of the Zariski closure of $Z_u$  in $V$.
\end{df}

We note that this definition differs from the one given in
\cite{CHKS} where we also included the critical points
in $\, V \backslash \mathcal{U}\,$ and we counted them
with multiplities.

Let $\{g_1,g_2,\ldots,g_r\}$ be a set of homogeneous polynomials
generating the ideal $P$. We consider the Jacobian matrix
augmented by a row of ones:
\begin{equation} \label{augmented}
J \quad = \quad \bmatrix
& 1 & 1 & \cdots & 1 \\
& {\partial g_1}/{\partial p_0} & {\partial g_1}/{\partial p_1} &
\cdots &
{\partial g_1}/{\partial p_n} \\
& {\partial g_2}/{\partial p_0} & {\partial g_2}/{\partial p_1} &
\cdots &
{\partial g_2}/{\partial p_n} \\
&  \vdots & \vdots & \ddots & \vdots & \\
& {\partial g_r}/{\partial p_0} & {\partial g_r}/{\partial p_1} &
\cdots & {\partial g_r}/{\partial p_n}
\endbmatrix .
\end{equation}
We multiply $J$ by the diagonal matrix whose entries are the
unknowns to get
$$ \tilde J \quad = \quad J \cdot {\rm diag}(p_0,p_1,\ldots,p_n) .$$

\begin{prop} \label{jc}
A point $p \in \mathcal{U}$ is in the likelihood locus $Z_u$ if and 
only if the 
data vector $(u_0, \dots, u_n)$ is in the image of 
the transpose matrix ${\tilde J}^T(p)$.
\end{prop}

\begin{proof}
Let $V_{\rm aff}$ be the affine subvariety of $\cc^{n+1}$ defined by 
$P + \langle \sum p_i - 1 \rangle$.  The Jacobian
of $V_{\rm aff}$ is the matrix $J$. The likelihood function
$L$ has no poles or zeros on $\mathcal{U}$, so the critical
points of $L$ are the same as the critical points of $\log(L) = \sum
_i u_i \log p_i$ on $V_{\rm aff}$. A point $p \in \mathcal{U}$ is a critical point of $\log(L)$ 
if and only if $ {\rm dlog} (L)(p) = (\frac{u_0}{p_0}, \dots,
\frac{u_n}{p_n})$ is in the image of $J^T(p)$. As $p_i \neq 0$
on $\mathcal{U}$, this is equivalent to $u = (u_0, \dots, u_n)$ being in the image
of ${\tilde J}^T(p)$.
\end{proof}

Our algorithm for computing the likelihood ideal $I_u$ will be derived from
Proposition \ref{jc}. First, however, let us show that
$I_u$ is always artinian for generic $u$.
Hence the colength of $I_u$ is constant
for almost all data $u$. This 
number is the {\em maximum likelihood (ML) degree} of the
projective variety $V$.

\begin{prop} \label{genfin} Let $\mathcal{I} \subset \mathcal{U} \times
\mathbb{P}^n$ be the incidence variety consisting of pairs $(p,u)$
where $p \in Z_u$. Then $\mathcal{I}$ is the projectivization of a
vector bundle over $\mathcal{U}$ and $\dim \mathcal{I} = n$. In particular, 
$Z_u$ is either empty or finite for generic $u$.
\end{prop}

\begin{proof} Let $c$ be the codimension of $V$. 
For every $p \in \mathcal{U}$ the
matrix $J(p)$ and hence the matrix ${\tilde J}(p)$ with their first 
rows removed  have rank $c$. Multiplying
${\tilde J}$ by the vector of ones yields $(\sum p_i, 0, \dots, 0)$. In
particular, for any $p \in \mathcal{U}$, the first row is linearly
independent of the remaining rows, and ${\tilde J}(p)$ has rank $c+1$. Thus
the set of all $u$ in the image of ${\tilde J}^T(p)$ is a vector space of
dimension $c+1$, and hence $\mathcal{I}$ is the projectivization
of a vector bundle of rank $c+1$ over $\mathcal{U}$. It follows that $\dim
\mathcal{I} = \dim \mathcal{U} + c = n$. Projecting onto the second factor, the
generic fiber must either be empty or of dimension $0$.
\end{proof}

\begin{ex} \label{circle}
Let $n = 2$ and 
$\, P = \langle p_0^2 + p_1^2 + p_2^2 - 2p_0p_1 - 2 p_0 p_2 - 2 p_1 p_2 \rangle$. 
The model $V$ is a circle in the triangle $\Delta_2$ which is
tangent to the three edges of $\partial \Delta_2$. 
The critical ideal $I_u \subset \rr[V]$ contains the cubic polynomial
\begin{equation}
\label{circlecubic}
u_2 p_0^2 p_1\,-\,
u_1 p_0^2 p_2 \,-\,
u_2 p_0 p_1^2 \,+\,
u_1 p_0 p_2^2\,+\,
u_0 p_1^2 p_2\,-\,
u_0 p_1 p_2^2.
\end{equation}
If $u_0,u_1,u_2$ are distinct, then this cubic curve intersects the circle
in six points, but only three of them lie in $\,\mathcal{U} $, which is
the part of the circle in the interior of the triangle $ \Delta_2 $.
The ML degree of the circle $V$ is three. Hence 
our problem (\ref{mleproblem}) can be solved in terms of radicals:
use {\em Cardano's formula} to express
each of the three points in $Z_u$  as
a function of the data $u_0,u_1,u_2$. \qed
\end{ex}

In Example \ref{circle}, the incidence variety $\mathcal{I}$ is the surface
in $\mathcal{U} \times \pp^2$ defined by 
$(\ref{circlecubic})$, which is regarded as a bihomogeneous
equation of degree $(3,1)$ in $(p,u)$. 

\section{Complete Intersections}

Here we consider the case when our model $V \subset \pp^n$ is a
complete intersection. This means that the codimension $c$
of $V$ coincides with the number $r$ of generators of the ideal $P$.
As before, we write $\, P \, = \,\langle g_1,g_2,\ldots,g_r \rangle $.
Let $d_i$ be the degree of the homogeneous polynomial $g_i$.
 Let $D$ denote the sum of all
monomials of degree at most $n-r$ in $r$ unknowns evaluated at
$d_1,d_2,\ldots,d_r$:
\begin{equation}
\label{thom} D \quad \,\, = \quad
\sum_{i_1 + i_2 + \cdots  + i_r \leq n-r} \!\! d_1^{i_1} d_2^{i_2}
\cdots d_r^{i_r}.
\end{equation}

\begin{thm} \label{CI} The ML degree of the model $\,V$ is bounded above by
$ D d_1 d_2 \cdots d_r $.
Equality holds when $\,V$ is a generic complete intersection, that is,
when the coefficients of the defining polynomials $g_1,g_2,\ldots,g_r$ are chosen at random.
\end{thm}

To illustrate this formula, let us consider some special cases.
First, suppose that our model $V$ is a hypersurface
 ($r=1$) defined by one homogeneous polynomial $g = g_1 $
  of degree $d = d_1$. Then the ML degree of $V$ is at most
  \begin{equation}
  \label{hypersurface} 
  d \cdot D \qquad = \qquad
 d \cdot \frac{d^n-1}{d-1}. 
 \end{equation}
In Example \ref{circle}, we considered the case of a quadric
in the plane $(d=n=2)$ having ML degree three.
The upper bound (\ref{hypersurface}) equals six, and 
this is indeed the ML degree of a general quadric.
Two special quadrics of statistical interest are the
{\em Hardy-Weinberg curve} $\, p_1^2 = 4 p_0 p_2\,$
and its  cousin $\,p_1^2 = p_0 p_2 $. The ML degrees of these 
two special models are one and two respectively.

Another noteworthy special case arises when
$V$ is a linear space of codimension $r$ in
$\pp^n$, i.e., $d_1 = \cdots = d_r = 1$.
Here the open set $\mathcal{U}$ is the (complexified)
complement of an arrangement of
$n+1$ hyperplanes in $\rr^{n-r}$, and
the ML degree equals the number of 
bounded regions of the (real) arrangement
\cite[\S 4]{CHKS}. If $V$ is generic then
the number of bounded regions equals
\begin{equation}
\label{linearspace}
d \cdot D  \quad = \quad 1 \cdot \binom{n}{r} \quad = \quad
\binom{n}{r}.
\end{equation}
An important statistical application of such linear models
is discussed in \cite{BR}.

\begin{proof}[Proof of Theorem \ref{CI}]
We first consider the case when the $g_i$ are generic forms
and $u$ is generic. By Bertini's Theorem, the generic complete
intersection $V$ is smooth. 
All critical points of the likelihood function $L$ on $V$
lie in the dense open subset $\mathcal{U}$, and the set $Z_u$ of critical points
is finite, by Proposition \ref{genfin}.

 Consider the following $(r+2) \times (n+1)$-matrix with entries in $\rr[p_0,\ldots,p_n]$:
 $$ \bmatrix u \\ \tilde J \endbmatrix \quad = \quad
\bmatrix
& u_0 & u_1 & \cdots & u_n \\
& p_0 & p_1 & \cdots & p_n \\
&
p_0 \frac{\partial g_1}{\partial p_0} &
p_1 \frac{\partial g_1}{\partial p_1} & \cdots &
p_n \frac{\partial g_1}{\partial p_n} \\
&
p_0 \frac{\partial g_2}{\partial p_0} &
p_1 \frac{\partial g_2}{\partial p_1} & \cdots &
p_n \frac{\partial g_2}{\partial p_n} \\
&  \vdots & \vdots & \ddots & \vdots & \\
&
p_0 \frac{\partial g_r}{\partial p_0} &
p_1 \frac{\partial g_r}{\partial p_1} & \cdots &
p_n \frac{\partial g_r}{\partial p_n}
\endbmatrix .
$$
 Let $W$ denote the determinantal variety in $\pp^n$ given by the vanishing of its
$(r+2) \times (r+2)$ minors. 
The codimension of $W$ is at most
$n-r$, which is a general upper bound for ideals of maximal
minors, and hence the dimension of $W$ is at most $r$.
Our genericity assumptions ensure that
the matrix $\tilde J(p)$ has maximal row rank $r+1$ for all $p \in V$.
Hence a point $p \in V$ lies in $W$
if and only if the vector $u$ is in the row span of $\tilde J(p)$.
Proposition \ref{jc} implies 
$$ Z_u \quad = \quad \mathcal{U} \,\cap \, W \quad = \quad V \, \cap \, W . $$
Since $Z_u$ is finite and $V$ has dimension 
$n-r$, we conclude that $W$ has the maximum possible codimension, 
namely $n-r$,  and that the
intersection of $V$ with the determinantal variety $W$ is transversal.
We note that $W$ is  Cohen-Macaulay, since $W$ has 
maximal codimension $n-r$, and 
ideals of minors  of generic matrices are Cohen-Macaulay.
B\'ezout's Theorem \cite[\S 8.4]{Ful} implies
$$
\hbox{ML degree} \quad = \quad {\rm degree}(V) \cdot {\rm
degree}(W) \quad = \quad d_1 \cdots d_r \cdot {\rm degree}(W). $$
The degree of the determinantal variety $W$  equals
the degree of the determinantal variety given by generic
forms of the same row degrees.   A special case of the
Thom-Porteous-Giambelli formula \cite[\S 14.4]{Ful} states that this degree is
the complete homogeneous symmetric function of degree $\,{\rm
codim}(W) = n-r\,$ evaluated at the row degrees of the matrix.
Here, the row degrees are $\,0,1,d_1, \ldots, d_r $, and the value
of that symmetric function is precisely $D$.
We conclude that ${\rm degree}(W) = D$.
This completes the proof that the ML degree
of the  generic complete intersection $\,V = \mathcal{V}(
g_1,\ldots,g_r ) \,$ equals
$\, D \cdot d_1 d_2 \cdots d_n $.

Suppose now that the $g_i$ are no longer generic. The ML degree of
$\,V = \mathcal{V}( g_1,\ldots,g_r )\,$ remains finite
by Proposition \ref{genfin}. The deformation argument in
 \cite[Theorem 22]{CHKS}
 implies that the ML degree of $V$ is at most
 $\, D \cdot d_1 d_2 \cdots d_n $.
\end{proof}

\section{Algorithms and Implementation}

We propose the following algorithm for deriving
the likelihood equations. 

\begin{alg} \label{implicit} (Computing the likelihood equations)  \rm \hfill \break
{\sl Input:} A homogeneous ideal $P \subset \rr[p_0,\ldots,p_n]$ and a vector
$u \in \nn^{n+1}$.\\
{\sl Output:} The likelihoood ideal $I_u$ of the model $V = \mathcal{V}(P)$ for the data $u$. 
\smallskip \\
{\sl Step 1:} Compute $c = {\rm codim}(V)$.
Let $Q$ be the ideal of
the singular locus of $V$, i.e.,  $Q$ is generated
by the $c \times c$ minors of the Jacobian matrix of $P$.
\\
{\sl Step 2:} Compute the kernel $M$ of the matrix 
$\,{\tilde J}\,$ over $\,\rr[V] = \rr[p_0,\! ..  ,p_n]/P$. 
{\sl Step 3:} Let $I'_u$ be the ideal in $\rr[V]$
generated by the polynomials $\, \sum_{i=0}^n u_i \cdot \phi_i $, where
the vectors
$\,(\phi_0,\ldots,\phi_n)\,$ 
run over a generating set of the module $M$.
\\
{\sl Step 4:}
The ideal $I_u$ equals the saturation 
$\bigl( I'_u : (p_0 \cdots p_n (\sum p_i) \cdot Q)^\infty \bigr)$.
\end{alg}

\begin{proof}[Proof of correctness]
By Proposition \ref{jc}, a point $p \in \mathcal{U}$ lies
 in $Z_u$ if and only if $u \cdot \phi(p) = 0$ 
for every $\phi(p)$ in the kernel of ${\tilde J}(p)$.
Since ${\tilde J}(p)$ has constant rank 
for all $p \in \mathcal{U}$, generators of the vector space 
${\rm kernel}_{\cc}({\tilde J}(p))$ are gotten by specializing
generators of the module
$\,M = {\rm kernel}_{\rr[V]}(\tilde J)$.
This shows that the ideal $I'_u$ vanishes on $Z_u$. 
Now, let $f$ be a polynomial in the saturation of Step 4, i.e.
$f \cdot (p_0 \cdots p_n \cdot g)^k \in I'_u$ 
for some $g \in Q$ and $k \in \mathbb{N}$.
Since this product vanishes on $Z_u$, 
the  polynomial $f$ vanishes on $Z_u$, and hence $f \in I_u$.

Conversely, for any $g \in Q$, the module 
$M$ has a free basis over the localization $\,\rr[V]_{g \cdot p_0 \cdots p_n}$.
Any element $f$ of $I_u$ is a linear combination of the dot product of $u$ 
with these free generators with coefficients in $\,\rr[V]_{g \cdot p_0 \cdots p_n}$.
By clearing denominators we get a polynomial which is 
a polynomial linear combination of the generators of $I'_u$.   
This shows that $f$ is in the saturation.
\end{proof}

\begin{rem} The ML degree of $V$ is computed by running
Algorithm \ref{implicit} for a generic vector $u \in \rr^{n+1}$.  We simply
output the colength of $I_u$ after Step 4. 
\end{rem}

A key feature of Algorithm \ref{implicit} is that Step 1 and Step 2 are
independent of the data $u$, so they need to be run only once
per model. Moreover, these preprocessing steps can be 
enhanced by doing the saturation of Step 4 already once at the
level of the module $M$, i.e., after Step 2 one can replace $M$ by
$$ \tilde M \quad := \quad \bigl( M : (p_0 \cdots p_n \cdot Q)^\infty \bigr) \quad = \quad
\rr[V]_{g \cdot p_0 \cdots p_n}\cdot M \,\,\,\cap \, \,\, \rr[V]^{n+1}. $$
For any particular data vector $u \in \nn^{n+1}$, one can then use
either $M$ or $\tilde M$ in Step 3 to define $I'_u$. The remaining
saturation in Step 4 requires some tricks in order to run efficiently.
We found that, for many models and most data, it suffices to
saturate only once with respect to a single polynomial, as follows:

\medskip 

\noindent {\sl Step 4'}: Pick a random $(c+1) \times (c+1)$-submatrix
of $\tilde J$ and let $h$ be its determinant. With some luck,
the likelihood ideal  $I_u$ will be equal to
$(I'_u : h)$.

\medskip

Here is one more useful variant.
When $V$ is a complete intersection, one can jump
directly to Step 3 and replace $I'_u$ by the determinantal variety
 $W$ in the proof of Theorem \ref{CI}. Thus, instead of $I'_u$
 we simply take the ideal of $(r+2) \times (r+2)$ minors of the matrix
$\, \bmatrix u \\ \tilde J \endbmatrix$. This variant is usually
slower than Algorithm \ref{implicit}, but it is sometimes 
faster when the codimension $r$ is small.

\smallskip

Here is one more comment concerning Step 2. Suppose our
computer algebra system does not support linear algebra
over quotient rings (such as $\rr[V]$). Then we can
implement Step 2 over the
polynomial ring $\rr[p_0,\ldots,p_n]$ as follows.
Instead of computing the kernel of
the $ (r+1)\times (n+1)$-matrix ${\tilde J}$, we compute the kernel of
the $ (r\!+ \!1)\times (n +1+r + r^2)$-matrix $[\,{\tilde J} \,|\, G\,]$, where
$$
G \quad = \quad \bmatrix
g_1 & \cdots & g_r &
 0  & \cdots &  0  & \cdots &
 0  & \cdots &  0  \\
 0  & \cdots &  0 &
g_1 & \cdots & g_r & \cdots &
 0  & \cdots &  0  \\
 \vdots & \vdots & \vdots &
 \vdots & \vdots & \vdots &  \ddots &
 \vdots & \vdots & \vdots \\
 0  & \cdots &  0 &
 0  & \cdots &  0 & \cdots &
g_1 & \cdots & g_r   \\
\endbmatrix
$$
Take the first $n+1$ coordinates from the generators
of  the kernel of $[\,{\tilde J} \,|\, G\,]$. These vectors 
generate the module $M $, so they
can be used in Step 3.

\medskip

Recall that our objective is to compute maximum likelihood estimates.

\begin{alg} \label{numerical} (Computing the local maxima of the likelihood function)  \rm \hfill \break
{\sl Input:}  The likelihood ideal $\,I_u$\, for the model $\,V\,$ and the data $\,u$. \\
{\sl Output:} The list of all local maxima for the optimization problem (\ref{mleproblem}).
\smallskip \\
{\sl Step 1:} If $\dim(I_u) = 0$ for the given data $u$, compute 
the solution set $Z_u$ numerically using Gr\"obner bases and eigenvalue methods,
as in \cite[\S 2]{CLO}.
\\
\indent
For each positive solution $\, p^* \in Z_u \cap V_{> 0} \, $ perform the following steps: \\
{\sl Step 2:}
 Solve the linear system ${\tilde J}^T(p^*) \cdot \lambda= u$ to get
Lagrange multipliers $\lambda^*_i$. 
The {\em Lagrangian}
 $\, {\mathcal L} := \log(L(p)) - \sum_{i=1}^r \lambda^*_i g_i (p)\,$
is a function of $p$.
\\ 
{\sl Step 3:} Compute the Hessian $H(p)$ of the Lagrangian ${\mathcal L}(p)$.
Compute the restriction of $H(p^*)$ to the tangent space  
$\mathrm{kernel}(J(p^*))$ of $V$ at the point $p^*$.
\\
{\sl Step 4:} If the restricted $H(p^*)$ in Step $3$ is negative definite, then output $p^*$
with its log-likelihood $\log(L(p^*))$ and
 the eigenvalues of the restricted $H(p^*)$.
\end{alg} 

We implemented Algorithms \ref{implicit} and \ref{numerical} in the computer algebra package 
{\sc Singular} \cite{GPS01}. 
The input is a homogeneous ideal $P$ in a polynomial ring and a data
vector $u$. The output is the ML degree and a list of all positive local maxima
$p^*$  and their certificates, namely the (negative) eigenvalues of the
Hessian $H(p^*)$.   {\sl Step 4} of Algorithm \ref{numerical} 
uses the well-known second order optimality conditions in 
nonlinear optimization, see for instance \cite{NW}.

All computational results
to be reported in the following sections were obtained
using this implementation.
 We also implemented Algorithm
\ref{implicit} in {\sc Macaulay 2} \ \cite{M2}.
This 
independently confirmed the reported ML degrees.

\section{Small Matrix Models}   

Determinantal varieties are natural objects both in algebraic geometry
and in statistics. In this section we discuss likelihood
equations, ML degree, and local maxima for some models 
specified by rank conditions on $3 \times 3$ matrices.

\begin{ex} \label{mixBern}
Consider the mixture model for Bernoulli random variables
discussed in the Introduction. This model is given by
 the determinant of
$$   {\bf P} \,\,\, = \,\,\,
\bmatrix
12 p_0 & 3  p_1 & 2 p_2 \\
 3 p_1 & 2 p_2 & 3 p_3 \\
 2 p_2 & 3 p_3 & 12 p_4
\endbmatrix.
$$
The ML degree of this model is twelve, and all twelve solutions
to the critical equations can be real.
In our experiments we found  that at most six of these solutions
are real and positive, and three of those can be local maxima. 
 A data vector for which the function  (\ref{likelihoodfct}) has
 three positive local maxima is 
$$ u \quad = \quad (u_0, u_1, u_2, u_3, u_4) \quad = \quad (51, 18, 73, 25, 75). $$ 
\end{ex}

\begin{ex} \label{SweetThirteen} 
Consider the general $3 \times 3$-matrix with indeterminate entries
$$  {\bf P} \,\,\, := \,\,\,
\bmatrix
p_{00} & p_{01} & p_{02} \\
p_{10} & p_{11} & p_{12} \\
p_{20} & p_{21} & p_{22}
\endbmatrix.
$$
The prime ideal of  $2 \times 2$ minors of this matrix represents two
independent ternary random variables. This model has ML degree
one. In other words, the critical equations have 
a unique (positive) solution for a given $3 \times 3$ data
matrix $U$. This maximum likelihood
estimate  is a rational function in $U$, namely, it is
the unique matrix of rank one with the same row and column sums as $U$.
This example is an instance of a  
{\em decomposable graphical model} and it is known that the ML degree of
such a  model is always one \cite{GMS}.

Continuing with our example, let $P$ be the principal ideal generated 
by the $3 \times 3$ determinant of ${\bf P}$. 
This is the mixture model for two pairs of
independent ternary random variables. The ML degree of this
mixture model equals $10$.
For a concrete numerical example consider the following data:
$$ U \,\,\, = \,\,\, \bmatrix 
16 & 17 & 7 \\
18 & 3 & 12 \\
1 & 8 & 16             \endbmatrix .$$
The likelihood ideal $I_U$ has four imaginary zeros and six real zeros,
all of which lie in the positive orthant. Three of these six
matrices are local maxima of the likelihood function. We list the 
three local maxima together with the values of the likelihood
function. The third matrix is the global maximum:
$$ \bmatrix 
.13887222 & .18906469 & .080226355 \\
.12570444 & .039074119 & .17195613 \\
.092566192& .057575479 & .10496037  \endbmatrix, \,\,\quad \log(L) = -207.0295890,
$$
$$
\bmatrix 
.14622787& .11633326& .14560213\\
.19982703& .046435565& .090472102\\
.011087957& .12294546& .12106862   \endbmatrix, \,\,\quad \log(L) = -202.9010713,
$$
$$
\bmatrix 
.20299213& .11762942& .087541717\\
.14331103& .096617294& .096806365\\
.010839697& .071467568& .17279478  \endbmatrix, \,\,\quad \log(L) = -202.6703908.
 $$ 

As was mentioned in the Introduction, the ML degree is very 
sensitive to even slight perturbations  to the ``natural'' coordinates
of the model. To illustrate this, let us scale the unknowns and
consider the new matrix
$$  {\bf P'} \,\,\, := \,\,\,
\bmatrix
\alpha_{00} p_{00} & \alpha_{01} p_{01} & \alpha_{02} p_{02} \\
\alpha_{10}  p_{10} & \alpha_{11}  p_{11} & \alpha_{12} p_{12} \\
\alpha_{20} p_{20} & \alpha_{21} p_{21} & \alpha_{22} p_{22}
\endbmatrix.
$$
where the $\alpha_{ij}$ are random real numbers. 
It turns out that the ML degree of the ideal of $2 \times 2$ minors of 
${\bf P'}$ jumps  to six. The ML degree of the ideal generated by the
determinant of ${\bf P'}$ jumps from $10$ to $39$ after this change. \qed
\end{ex} 

\begin{ex} \label{gen3x3}
Consider the following symmetric $3 \times 3$-matrix:
$$  {\bf P} \,\,\, := \,\,\,
\bmatrix
2p_{00} & p_{01} & p_{02} \\
p_{01} & 2p_{11} & p_{12} \\
p_{02} & p_{12} & 2p_{22}
\endbmatrix.
$$
The ideal of $2 \times 2$ minors of this matrix represents two
independent identically distributed ternary random variables. This
model has ML degree $1$. Again, the mixture model for two copies
of the previous model is specified by the determinant of ${\bf P}$. The
ML degree of this mixture model equals $6$.

Note how these ML degrees change if we replace ${\bf P}$
by the scaled  matrix
$$  \qquad {\bf P'} \,\,\, := \,\,\,
\bmatrix
\alpha_{00} p_{00} & \alpha_{01} p_{01} & \alpha_{02} p_{02} \\
\alpha_{01} p_{01} & \alpha_{11} p_{11} & \alpha_{12} p_{12} \\
\alpha_{02} p_{02} & \alpha_{12} p_{12} & \alpha_{22} p_{22}
\endbmatrix
\qquad \hbox{(with $\alpha_{ij}$ random reals)}.
$$
The ideal of $2 \times 2$ minors of ${\bf P'}$ has ML degree $4$, which
is the degree of the corresponding Veronese surface. The first
secant variety of the Veronese surface is given by the determinant
of ${\bf P'}$. That model has ML degree $16$. \qed
\end{ex}

We close this section with a table of ML degrees for seven
determinantal varieties. The first and second columns are
Examples \ref{SweetThirteen} and \ref{gen3x3} respectively.
The first row indicates the original ideal
in the statistically natural coordinates. The second row 
refers to the (``scaled'') ideal gotten from the first (``unscaled'') ideal by
 generically scaling the coordinates. The column
 ``$ 3  \times 4$'' refers to the maximal minors of a $3 \times 4$-matrix,
 ``$ 3 \!  \times \! 3 coin$'' is Example \ref{mixBern}, and
 ``$ 3 \! \times \! 4 coin$'' is a similar problem where the coin
is tossed five times in a row instead of four, and $P$ is the ideal of
the $3 \times 3$ minors of the matrix 
$$
\bmatrix
10p_0 & 2p_1&  p_2 &   p_3 \\
2p_1 &  p_2 &  p_3 & 2p_4 \\
p_2 &  p_3 & 2p_4 & 10p_5 
\endbmatrix.
$$
Finally, $G(m,n)$ 
is
 the Pl\"ucker ideal of the Grassmannian
of $m$-planes in $\cc^n$.

$$\begin{array}{|c||c|c|c|c|c|c|c|}
\hline
{\rm Model} & 3 \times 3  & 3 \! \times \! 3sym & 3 \times 4 & 3 \! \times 
\! 3coin  & 3 \! \times \! 4coin & G(2,4) & G(2,5) \\
\hline
{\rm unscaled} &  10    &            6      &  26         &  12                      
&            39       &   4     &  22 \\
\hline
{\rm scaled} &      39    &           16     &  164     &  16                       
&    54                 &   6   &   52   \\
\hline
\end{array}
$$

\medskip

\noindent
We close with two open problems, aimed at experts in enumerative geometry.

\begin{pr} Find an explanation for all the ML degrees stated above.  
\end{pr}

\begin{pr} Characterize all models whose ML degree 
is one.
\end{pr}

\section{Jukes-Cantor Models in Phylogenetics}

The study of ``analytic solutions''  for maximum likelihood 
estimation has a long tradition in phylogenetics \cite{Fel},
where one considers evolution models for DNA sequence
data, and maximum likelihood is used to find the best 
phylogenetic tree that explains the evolution of the taxa
under consideration. Maximum likelihood  
is also used to estimate the branch lengths of the
reconstructed trees. 
Here we examine the widely used Jukes-Cantor models,
with emphasis on the cases studied by
Chor et al. \cite{CHS, CHPH, CKS}
and Sainudiin \cite{Sai}.

We  use the notation of Sturmfels and Sullivant \cite{StSu},
first for binary data and later (in Example \ref{jcdna})  
for DNA data. Let us start out with Example 3 in \cite{StSu}.
We consider any tree $T$ with three leaves and the Jukes-Cantor
model with unknown root distribution.  This is
equivalent to considering trees with four leaves and uniform root
distribution. Each tree topology $T$ specifies a model for three binary
random variables. The joint probabilities are represented by
unknowns $p_{ijk}$, for $i,j,k  \in \{0,1\}$. The data is given
as a $2 \times 2 \times 2$-table $u = (u_{ijk})$ whose entries
record the number of occurrences of any particular column pattern
among three aligned binary sequences.
We perform the linear change of coordinates given by the {\em Fourier transform}:
\begin{eqnarray*}
& q_{000} \quad = \quad
 p_{000}+p_{001}+p_{010}+p_{011}+p_{100}+p_{101}+p_{110}+p_{111} , \\
& q_{001} \quad = \quad p_{000}-p_{001}+p_{010}-p_{011}+p_{100}-p_{101}+p_{110}-p_{111} , \\
& q_{010} \quad = \quad p_{000}+p_{001}-p_{010}-p_{011}+p_{100}+p_{101}-p_{110}-p_{111} , \\
& q_{011} \quad = \quad p_{000}-p_{001}-p_{010}+p_{011}+p_{100}-p_{101}-p_{110}+p_{111} , \\
& q_{100} \quad = \quad p_{000}+p_{001}+p_{010}+p_{011}-p_{100}-p_{101}-p_{110}-p_{111} , \\
& q_{101} \quad = \quad p_{000}-p_{001}+p_{010}-p_{011}-p_{100}+p_{101}-p_{110}+p_{111} , \\
& q_{110} \quad = \quad p_{000}+p_{001}-p_{010}-p_{011}-p_{100}-p_{101}+p_{110}+p_{111} , \\
& q_{111} \quad = \quad
p_{000}-p_{001}-p_{010}+p_{011}-p_{100}+p_{101}+p_{110}-p_{111} .
\end{eqnarray*} The advantage of this transformation is that the
defining ideal $P$ of any Jukes-Cantor model becomes a {\em toric
ideal } in the {\em Fourier coordinates} $q_{ijk}$.

\begin{ex} \label{claw}
Let $T = K_{1,3}$ be the claw tree with three edges
attached to the root. Then our  model is a complete intersection of
codimension $3$ in $\pp^7$:
$$ P  \quad = \quad \langle \,
q_{001} q_{110} - q_{000} q_{111}, \, q_{010} q_{101} - q_{000}
q_{111},\, q_{100} q_{011} - q_{000} q_{111} \,\rangle.
$$
Our problem is to solve the following constrained optimization
problem:
$$ {\rm maximize}\, \prod_{i=0}^1 \prod_{j=0}^1 \prod_{k=0}^1
p_{ijk}^{u_{ijk}}  \quad \hbox{subject to} \,\, p \,\in \mathcal{V}(P)
  \,\, \hbox{and } \,\,\sum_{ijk} p_{ijk} = 1. $$
Algorithm \ref{implicit} easily derives the likelihood equations,
and it reports that, for random data $u$, the equations have  $92$
distinct complex solutions. In short, the Jukes-Cantor binary
model on the claw tree $K_{1,3}$ has ML degree $92$.
 \qed
\end{ex}

\begin{ex} \label{nonclaw}
Suppose that $T$ is one of the three trivalent trees, for
instance, the one where  the leaves $1$ and $2$ are split from the
leaf $3$. This model is a complete intersection of codimension two.
The ideal of model invariants is
$$ P \quad = \quad \langle\,\,
q_{001} q_{110} - q_{000} q_{111},\, q_{010} q_{101} - q_{100}
q_{011} \,\rangle. $$
This model has dimension $5$ and ML degree $14$. 
We found many instances where two of the
$14$ complex solutions to the likelihood equations are
local maxima in the probability simplex $\Delta_7$,
thus confirming the results of \cite{CHPH}.

The authors of \cite{CKS} studied the
three-dimensional submodels gotten
by assuming the {\em molecular clock}
hypothesis. There are two combinatorial types:
\begin{eqnarray*} P_{fork} \quad = &
 \langle \, q_{100} - q_{101},
              \, q_{011} - q_{101}, 
             \, q_{010} - q_{101} , \, q_{001} q_{110} - q_{000} q_{111} \,\rangle \\
P_{comb} \quad = &       \langle \,                   
             q_{010} - q_{100},\,
              q_{001} - q_{100},\, q_{011} - q_{101},\, q_{100} q_{110} - q_{000} q_{111} \,\rangle.
\end{eqnarray*}
The ideal $P_{fork}$ has ML degree one, and the ideal
$P_{comb}$ has ML degree nine.
It was shown in \cite{CKS} that the local maximum in
$\Delta_7$ is unique for $P_{comb}$.
\qed
\end{ex}

Each rooted tree with leaves $\{1,2,3\}$ is specified by its {\em
split system}, which is a collection $\Sigma$ of splits of the set
$\{0,1,2,3\}$ into two non-empty parts. Here $0$ represents the
root. The number of splits equals the dimension of the model. The
split systems representing the trees in  Example \ref{claw} and
Example \ref{nonclaw} are
\begin{eqnarray*}
& \Sigma_{\rm \ref{claw}} \,\, = \,\,\bigl\{ \{0,123\},
\{1,023\}, \{2,013\} ,\{3,012\} \bigr\}, \qquad
\qquad \\
& \Sigma_{\rm \ref{nonclaw}}\,\, = \,\,\bigl\{\{03,12\},
\{0,123\}, \{1,023\}, \{2,013\} ,\{3,012\} \bigr\}.
\end{eqnarray*}
David Bryant \cite{Bry} proposed to generalize phylogenetic models from trees
to arbitrary splits graphs. Jukes-Cantor models for splits graphs
are likely to become important for applications. Here is the
simplest non-tree example:

\begin{ex} \label{splitmodel}
We add one more split to $\,\Sigma_{\rm \ref{nonclaw}}\,$ 
to get the split system
$$ \Sigma_{\rm \ref{splitmodel}}
\,\,\, = \,\,\,
\bigl\{\{01,23\},\{03,12\}, \{0,123\}, \{1,023\}, \{2,013\}
,\{3,012\} \bigr\}. $$ 
The resulting Jukes-Cantor model is a hypersurface of degree four in $\pp^7$:
$$ P \quad = \quad \langle \,q_{000} q_{010} q_{101} q_{111} -
q_{001} q_{110} q_{110} q_{011}\, \rangle. $$ If we rewrite this
quartic in terms of the probabilities $p_{ijk}$ then we get a
polynomial with $40$ terms. The ML degree of this model equals
$326$.   Note that this is still a lot smaller than the upper bound
 $21,844$ given by (\ref{hypersurface}).
\qed
\end{ex}

All of the phylogenetic models whose likelihood equations have
been analyzed so far assumed binary characters. For applications
in biology, models on four character states (A, C, G and T) are
more important. We next present a detailed analysis of the
smallest non-trivial Jukes-Cantor DNA model.

\begin{ex} \label{jcdna}
Consider the Jukes-Cantor DNA model on a tree with three leaves
and uniform root distribution. The number of 
observable states is $4^3 = 64$
but it turns out that there are only five distinct probabilities.

We may  assume that the tree is the claw tree $ K_{1,3}$. The
model parameters $\pi_1,\pi_2,\pi_3$ are the probabilities of
changing from any letter (A,C,G or T) to any other letter when
passing from the root to the leaves $1,2,3$. We write $\,\theta_i
= 1 - 3  \pi_i$ for the probability of not changing the letter.
Let $p_{123}$ be the probability of observing the same letter at
all three leaves, $p_{ij}$ the probability of observing the  same
letter at all leaves $i,j$ and a different one at the third leaf,
and $p_{dis}$ the probability of seeing three distinct letters.
Then
\begin{eqnarray*}
 p_{123} \quad  = & \quad
\theta_1 \theta_2 \theta_3 \,\, + \,\,
3 \pi_1 \pi_2 \pi_3 , \\
 p_{dis} \quad   = & \quad
6 \theta_1 \pi_2 \pi_3 \,+\, 6 \pi_1 \theta_2 \pi_3 \,+\, 6 \pi_1
\pi_2 \theta_3 \,+\,
6 \pi_1 \pi_2 \pi_3 , \\
 p_{12} \quad  = & \quad
3 \theta_1 \theta_2 \pi_3 \,+ \, 3 \pi_1 \pi_2 \theta_3 \,+ \,
6 \pi_1 \pi_2 \pi_3 , \\
 p_{13} \quad  = & \quad
3 \theta_1 \pi_2 \theta_3  \,+ \, 3 \pi_1 \theta_2 \pi_3 \,+ \,
6 \pi_1 \pi_2 \pi_3 , \\
 p_{23} \quad  = & \quad
3 \pi_1 \theta_2 \theta_3 \,+ \, 3 \theta_1 \pi_2 \pi_3 \,+ \, 6
\pi_1 \pi_2 \pi_3 .
\end{eqnarray*}
Here we can either set $\theta_i = 1 - 3\pi_i$, or we  can also
regard $(\theta_i : \pi_i)$ as homogeneous coordinates for
$\pp^1$. The above formulas define a
 map $\pp^1 \times \pp^1 \times \pp^1 \rightarrow \pp^4$,
and our model $V$ is the image of this map. Its defining ideal
equals
$$ P \quad = \quad \langle \, q_{000} q_{111}^2 \, - \, q_{011} q_{101} q_{110}
\, \rangle . $$ Here the $q_{ijk}$ are the Fourier coordinates
which are specified by
\begin{eqnarray*}
& q_{111} \,\, = \,\,
  p_{123} + \frac{1}{3} p_{dis} - \frac{1}{3} p_{12}
  - \frac{1}{3} p_{13} - \frac{1}{3} p_{23}  \,= \, (\theta_1 -
\pi_1)(\theta_2 - \pi_2)(\theta_3 - \pi_3)
\\
& q_{110} \,\, = \,\,
 p_{123} - \frac{1}{3} p_{dis} + p_{12} - \frac{1}{3} p_{13}
- \frac{1}{3} p_{23}
  \,\, = \,\,
(\theta_1 - \pi_1)(\theta_2 - \pi_2)(\theta_3 + 3 \pi_3)
\\
& q_{101}\,\, = \,\,
 p_{123} - \frac{1}{3} p_{dis} - \frac{1}{3} p_{12}
+ p_{13} - \frac{1}{3} p_{23} \,\, = \,\, (\theta_1 -
\pi_1)(\theta_2 + 3 \pi_2)(\theta_3 - \pi_3)
\\
& q_{011} \,\, = \,\, p_{123}  - \frac{1}{3} p_{dis} - \frac{1}{3}
p_{12}
 - \frac{1}{3} p_{13} + p_{23}
  \,\, = \,\,
(\theta_1 + 3 \pi_1)(\theta_2 - \pi_2)(\theta_3 - \pi_3)
\\
& q_{000} \,\,\, = \,\,\, p_{123}  + p_{dis} +  p_{12} +  p_{13} +
p_{23}
  \,\,\, = \,\,\,
(\theta_1 + 3 \pi_1)(\theta_2 + 3 \pi_2)(\theta_3 + 3 \pi_3)
\end{eqnarray*}
Algorithm \ref{implicit} reveals that the ML degree of this model
equals $23$. Using Algorithm \ref{numerical}  we were able to
confirm the global maximum reported in \cite[Section 5.2]{Sai}
on DNA sequence data for Chimpanzee, Gorilla, and Orangutan. The data used
in this example is
$$ (u_{123},\, u_{dis}, \, u_{12}, \, u_{13}, \, u_{23}) \quad = \quad
(700,\, 7, \, 100, \, 42, \, 46),$$
where there is a second local maximum present. Out of the $23$ solutions
to the critical equations $17$ are real, and $7$ are positive.
Our experiments show that there are data for which as many as four
positive local maxima exist.

The authors of \cite{CHS} study the two-dimensional submodel
gotten by assuming the {\em molecular clock}
hypothesis. 
This is the surface in $\pp^4$ defined by
$$ P_{clock} \quad = \quad 
\langle \,
q_{011} \, - \, q_{101} \, , \,\,
q_{000} q_{111}^2 \, - \, q_{011} q_{101} q_{110}
\, \rangle . $$
The ML degree of $P_{clock}$ is $11$, confirming the 
{\sc maple} computation in \cite{CHS}.
\qed
\end{ex}

\section{Likelihood Equations from Parametrization}

Consider a statistical model which is given parametrically as
the image of a polynomial map $\, f : \rr^d \rightarrow
\rr^{n+1} $. Each coordinate $f_i$ of $f$ is a polynomial in
model parameters $\theta = (\theta_1,\ldots,\theta_d)$, and we
have $\,f_0 + f_1 + \cdots + f_n = 1$.
This is usually the natural presentation coming from statistics,
and it is the setting  of \cite{CHKS}.
The parametric version of (\ref{mleproblem}) is the following
optimization problem:
\begin{equation}
\label{paramax}  {\rm Maximize} \, \,
f_0(\theta)^{u_0} f_1(\theta)^{u_1} \cdots f_n(\theta)^{u_n}, 
\end{equation}
where $u  = (u_0,\ldots,u_n)$ is a vector of positive integers and
$\theta$ runs over an open subset of $\rr^d$. The {\em critical
equations} for this optimization problem are
\begin{equation}
\label{critical} \sum_{i=0}^n \frac{u_i}{f_i} \frac{\partial
f_i}{\partial \theta_j} \,\,\, = \,\, 0 \qquad \hbox{for}\,\, j =
1,\ldots,d .
\end{equation}

In this section we show how to solve these equations directly.
In our experience, the algorithms of Section 4 are generally
preferable if the ideal $P$ of algebraic relations among the $f_i$ is known.
But  sometimes the parametric algorithm presented below is
quite useful as well. Theorem \ref{thesame} 
says that (under reasonable assumptions)
both methods  produce the same answer.

 We consider the Zariski
open set $\mathcal{U}_f$ in $\mathbb{C}^d$ where none of the
$f_i$ are zero. The critical locus $Z_u$ is defined in $\mathcal{U}_f$ by the
vanishing of the equations \eqref{critical}. Let $J_u$ be the ideal in
$\rr[\theta] = \rr[\theta_1, \ldots , \theta_d]$ whose variety is the Zariski
closure of $Z_u$ in all of $\mathbb{C}^d$. We call $J_u$ the {\em
parametric likelihood ideal} for (\ref{paramax}).

Of course we can obtain $J_u$ by computing the ideal of the numerators
of the equations \eqref{critical} and then saturating by the product
of the $f_i$. This has the disadvantages of being quite slow in
practice and requiring a separate computation for each choice of
$u$. We propose the following method instead.

\begin{alg} \label{para}  (Computing the parametric likelihood equations)  \rm \hfill \break
{\sl Input:} Polynomials $f_0,\ldots,f_n \in \rr[\theta]$
with $\sum_i f_i = 1$ and a vector $u \in \nn^{n+1}$. \\
{\sl Output:} Generators of the parametric likelihood ideal $J_u \subset \rr[\theta]$. \smallskip \\
{\sl Step 1:} Compute generators for the kernel over $\rr[\theta]$ of the
matrix
\begin{equation}
\label{modelmatrix}
 M \quad = \quad \bmatrix
f_0 & 0 & \cdots & 0 & \frac{\partial f_0}{\partial \theta_1} & \cdots & \frac{\partial f_0}{\partial \theta_d} \\
0 & f_1 & \cdots & 0 & \frac{\partial f_1}{\partial \theta_1} & \cdots & \frac{\partial f_1}{\partial \theta_d} \\
  &     & \ddots &   &                                        & \ddots & \\
0 & 0   & \cdots &f_n&  \frac{\partial f_n}{\partial \theta_1} &
\cdots & \frac{\partial f_n}{\partial \theta_d}
\endbmatrix .
\end{equation}
{\sl Step 2:} For each generator $(\psi_0,
\dots, \psi_n, \xi_1, \dots, \xi_d)^T $  of
${\rm kernel}_{\rr[\theta]}(M) \,$
form the polynomial $\,\sum_{i=0}^n u_i \psi_i $. Let
$J'_u$ be the ideal generated by these polynomials. \smallskip \\
{\sl Step 3:}
The desired ideal is equal to the saturation 
$\,J_u  =  (J'_u : (f_0 f_1 \cdots f_n)^{\infty})$.
\end{alg}

The proof of correctness for Algorithm \ref{para} is straightforward
using the setup of \cite{CHKS}.
The kernel of the matrix $M$ is the module of {\em logarithmic
vector fields} along the hypersurface in $\cc^{d}$ defined by $f_0 f_1 \cdots
f_n$. It was shown in \cite[\S 7]{CHKS} that  $J'_u = J_u$ holds
under certain geometric hypotheses (namely, the map
$f$ factors through a smooth variety on which the $f_i$ represent
global normal crossing divisors). In general, we may still have to
saturate by  $\prod_i f_i$, but the generators of $J'_u$ 
are much closer to the ideal $J_u$ than the numerators of 
(\ref{critical}).

Unlike the implicit setting of Section 4, the ideal
$J_u$ need not be artinian even if $u$ is generic.
There can be positive-dimensional components of critical points
at the locus in $\theta$-space where the
parameterization fails to be smooth.

\begin{ex} \label{Bernoulli1} Let $d=3, n=4$ and consider the
example in the Introduction:
$$
f_0 \, = \, \pi  (1-s)^4 \, \, + \,\, (1-\pi)  (1-t)^4  \,,\,\ldots \,,\,
f_4 \, = \, \pi  s^4 \, \, + \,\,  (1-\pi)  t^4. $$
The kernel of the $5 \times 8$-matrix $M$ in (\ref{modelmatrix})
is minimally generated by $27$ vectors in $\rr[s,t,\pi]^8 $. 
We compute the parametric likelihood ideal
$J_u$ for generic $u$ using Steps 2 and 3 of Algorithm~\ref{para}.
It turns out that $J_u$ is not artinian and it has four associated
primes. The first is a
one-dimensional component:
$$ \langle \,
s - \hat U,\,  t-\hat U \,\rangle, \qquad \hbox{where} \,\quad
\hat U \,= \,\frac{u_1 + 2 u_2 + 3 u_3 + 4 u_4} {4
(u_0+u_1+u_2+u_3 + u_4 )} .$$ This component does not depend on
$\pi$ at all: This is the unique solution of the maximum
likelihood problem for the unmixed Bernoulli random variable. 
Next there are two components each of which contributes three
critical points:
\begin{eqnarray*}
& \langle \, \pi-1, \,s - \hat U, \,
\alpha_3 t^3 + \alpha_2 t^2 + \alpha_1 t + \alpha_0 \,\rangle  \\
\hbox{and} & \langle \,\,\, \pi\,,\,\, t - \hat U \,, \,\,
\alpha_3 s^3 + \alpha_2 s^2 + \alpha_1 s + \alpha_0 \, \rangle ,
\end{eqnarray*}
where the $\alpha_i$ are certain rational expressions in the
$u_j$. These critical points are extraneous. They 
 can be explained by noticing that the parameterization is singular when
 either $s=t$
or the mixing parameter $\pi$ equals $0$ or $1$.

 After saturating out these three extraneous components we
are left with an ideal $K_u$ which is prime over $\qq(u)$. It is
artinian and has $24$ complex zeros. These critical points come in
pairs $(\pi,s,t)$ and $(1-\pi,t,s)$. Removing this extra symmetry
confirms that the true ML degree of this model is $12$. \qed
 \end{ex}

This example suggests that we add one more step to Algorithm \ref{para}:

\smallskip

\noindent
{\sl Step 4: }
Let $Q$ be the ideal generated by the $d \times d$ minors of
the $(n+1) \times d$ 
Jacobian matrix $\,Df = \bigl(\partial f_i /\partial \theta_j \bigr)$.
Compute and output the saturation
\begin{equation}
\label{removesing}
 K_u \quad := \quad (\,J_u \, : \, Q^\infty ). 
\end{equation}

\smallskip

The variety $\mathcal{V}(Q)$ is the singular locus of the map $f$,
and the saturation (\ref{removesing}) removes all
components of the ideal $J_u$ that lie in this singular locus.
We close by relating the ideal $K_u$ to the ideal $I_u$
from Sections 2--4.

\begin{thm} \label{thesame} Let 
$f: \rr^d \rightarrow \rr^{n+1}$ be a polynomial
map whose image is defined by a homogeneous ideal $P$
as in Section 2. Suppose that  $f$ is  generically finite
of degree $\delta$, and
 the image of $\,\mathcal{U}_f \backslash
\mathcal{V}(Q)$ lies in the smooth locus of 
$V = \mathcal{V}(P)$. For generic $u \in \nn^{n+1}$, the variety
$\mathcal{V}(K_u)$ equals the preimage
of $\mathcal{V}(I_u)$.
In particular, $K_u$ is
artinian and its colength is $\delta$ times
the ML degree of $V$.
\end{thm}

\begin{proof} Let $g = (g_1, g_2, \dots, g_r)$ be the 
generators of $P$. Then we have $g \circ f = 0$. The Chain Rule implies
$\,J \cdot Df = 0$, where $J$ is the Jacobian of
$V_{\rm aff}$ as in (\ref{augmented}).
 The smoothness hypotheses guarantee that the rank of $J$
is $n+1-d$, while the rank of $Df$ is $d$, for all points
$p = f(\theta)$ where $\theta \in \mathcal{U}_f \backslash
\mathcal{V}(Q)$.
The dimension count shows that the image of $J^T$
equals the kernel of $Df^T$. 
More precisely, a vector $u$ lies in the kernel of $Df(\theta)^T$ if
and only if it lies in the image of $J(p)^T$ with 
$p = f(\theta)$. In light of Propositions
\ref{jc} and \ref{genfin},
this implies that, for $u$ generic,
every point $p \in \mathcal{V}(I_u)$ pulls back to 
$\delta$ points $\theta \in \mathcal{V}(K_u)$.
\end{proof}

\medskip

\noindent {\bf Acknowledgements:}
We are grateful to the Park City Mathematics Institute (PCMI,
July 2004) for providing us with the opportunity 
to work on  this project in the mountains of Utah.
Amit Khetan was supported by an NSF postdoctoral fellowship (DMS-0303292).
Bernd Sturmfels was supported by the Clay Mathematics Institute
and in part by the NSF (DMS-0200729).

\bigskip

\noindent {\bf Authors' addresses:}

\bigskip

\noindent  Serkan Ho\c sten, Department of Mathematics,
San Francisco State University, San Francisco, CA 94132, USA, \
{\tt serkan@math.sfsu.edu}

\medskip

\noindent Amit Khetan, Department of Mathematics,
University of Massachusetts, Amherst, MA 01002, USA, \
{\tt khetan@math.umass.edu}

\medskip

\noindent Bernd Sturmfels, Department of Mathematics,
University of California, \break Berkeley, CA 94720, USA, \
{\tt bernd@math.berkeley.edu}

\end{document}